\documentclass[a4paper,10pt]{amsart}

\usepackage{CJK}
\usepackage{amsmath}
\usepackage{amsthm,amssymb,latexsym}
\usepackage{amsmath,calligra}
\usepackage{url,ae}
\usepackage{t1enc}
\usepackage{eufrak}
\usepackage{amsfonts}
\usepackage{graphicx}
\usepackage{hyperref}
\usepackage{amscd}
\usepackage{mathrsfs}
\usepackage{bbm}
\usepackage{color}
\usepackage{txfonts}
\usepackage[all]{xy}

\usepackage{multirow,makecell}
\usepackage[format=hang,font=small,textfont=it]{caption}

\newtheorem{theorem}{Theorem}[section]
\newtheorem{corollary}[theorem]{Corollary}
\newtheorem{lemma}[theorem]{Lemma}
\newtheorem{proposition}[theorem]{Proposition}

\newtheorem{question}{Question}[section]

\theoremstyle{definition}
\newtheorem{definition}[theorem]{Definition}
\theoremstyle{remark}
\newtheorem{remark}[theorem]{Remark}
\theoremstyle{example}
\newtheorem{example}[theorem]{Example}

\newcommand{\be}{\begin{equation}}
\newcommand{\ee}{\end{equation}}
\newcommand{\bea}{\begin{eqnarray}}
\newcommand{\eea}{\end{eqnarray}}
\newcommand{\ben}{\begin{eqnarray*}}
\newcommand{\een}{\end{eqnarray*}}
\newcommand{\bet}{\begin{equation}
\begin{split}}
\newcommand{\eet}{\end{split}
\end{equation}}

\DeclareMathOperator{\Ext}{\mathscr{E}\text{\kern -3pt {\calligra\Large xt}}\,\,}

\begin{document}

\title[On Skoda's theorem for Nadel-Lebesgue multiplier ideal sheaves]
      {On Skoda's theorem for Nadel-Lebesgue multiplier ideal sheaves on singular complex spaces and regularity of weak K\"ahler-Einstein metrics}

\author{Zhenqian Li}

\date{\today}
\subjclass[2010]{14F18, 32L20, 32Q20, 32S05, 32U05, 32W05}
\thanks{\emph{Key words}. Multiplier ideal sheaves, plurisubharmonic functions, vanishing theorems, $\overline\partial$-equations, Skoda's $L^2$ division theorem, K\"ahler-Einstein metrics}
\thanks{E-mail: lizhenqian@amss.ac.cn}

\begin{abstract}
In this article, we will characterize regular points respectively by the local vanishing, positivity of the Ricci curvature and $L^2$-solvability of the $\overline\partial$-equation together with Skoda's theorem for Nadel-Lebesgue multiplier ideal sheaves associated to plurisubharmonic (psh) functions on any (reduced) complex space of pure dimension. As a by-product, we show that any weak K\"ahler-Einstein metric on \emph{singular} $\mathbb{Q}$-Fano/Calabi-Yau/general type varieties cannot be smooth, and that in general there exists no \emph{singular} normal K\"ahler complex space such that the K\"ahler metric is K\"ahler-Einstein on the regular locus.
\end{abstract}

\maketitle

\section{Introduction}

Throughout this note, all complex spaces are always assumed to be reduced and paracompact unless otherwise mentioned; we mainly refer to \cite{GR84, Richberg68} for basic references on the theory of complex spaces.

\subsection{Local vanishing for multiplier ideals}

The local vanishing theorem for the higher direct images of sheaves computing multiplier ideals plays an important role in complex geometry and algebraic geometry, by which many local/global properties of multiplier ideal could be deduced, e.g., the restriction theorem, Nadel vanishing theorem and Skoda's theorem for multiplier ideals and so on (cf. \cite{De10, La04}, etc.).

Let $(X,\omega)$ be a Hermitian complex space of pure dimension $n$ and $\varphi\in\text{QPsh}(X)$ be a quasi-psh function on $X$. Then we can define the Nadel-Lebesgue multiplier ideal sheaf $\mathscr{I}_\text{NL}(\varphi)$ associated to $\varphi$ on $X$ by the integrability with respect to the Lebesgue measure $dV_\omega$ (see Definition \ref{MIS_NL}), which coincides with the usual multiplier ideal sheaf $\mathscr{I}(\varphi)$ introduced by Nadel whenever $X$ is smooth. Let $\pi:\widetilde X\to X$ be any log resolution of the Jacobian ideal $\mathcal{J}ac_X$ of $X$, then it follows that the Nadel-Lebesgue multiplier ideal sheaf
$$\mathscr{I}_\text{NL}(\varphi)=\pi_*\left(\mathcal{O}_{\widetilde X}(\widehat K_{\widetilde X/X})\otimes\mathscr{I}(\varphi\circ\pi)\right).$$
When $X$ is smooth, the Mather discrepancy divisor $\widehat K_{\widetilde X/X}$ is nothing but the relative canonical divisor $K_{\widetilde X/X}:=K_{\widetilde X}-\pi^*K_{X}$ of $\widetilde X$ over $X$. Then, we have the following local vanishing for multiplier ideals (cf. \cite{La04, Matsu_morphism})
$$R^q\pi_*\Big(\mathcal{O}_{\widetilde X}(K_{\widetilde X/X})\otimes\mathscr{I}(\varphi\circ\pi)\Big)=0,\ \forall q\geq1.$$

Therefore, it is natural to ask whether we could establish a similar local vanishing result in the singular setting, i.e.,
$$R^q\pi_*\Big(\mathcal{O}_{\widetilde X}(\widehat K_{\widetilde X/X})\otimes\mathscr{I}(\varphi\circ\pi)\Big)=0,\ \forall q\geq1.$$
In the present note, one of our goals is to study the local vanishing in the context of Nadel-Lebesgue multiplier ideals. In particular, based on Skoda's division for Nadel-Lebesgue multiplier ideals, we will prove that such a local vanishing for Nadel-Lebesgue multiplier ideals is in fact equivalent to smoothness of the ambient space in some sense; see Theorem \ref{vanishing_regular} for a detailed statement.

\subsection{Skoda's ideal generation by $L^2$ estimates for the $\overline\partial$-equation}

In the classical works \cite{Skoda72a, Skoda78}, relying on the $L^2$ methods due to \cite{AV65, Hormander65, Hormander_book} in several complex variables, Skoda established an analytic criterion on the ideal generation by a given collection of holomorphic functions or sections. In the original proof of Skoda's ideal generation, as well as standard techniques in functional analysis for the argument on a priori estimate and solving $\overline\partial$-equation with $L^2$ estimates, he also developed special analytic techniques by restricting the domain of the $\overline\partial$-operator to an appropriate subspace of the usual $L^2$ space and inducing an $L^2$ estimate on this new operator.

As applications, Skoda's theorem is a crucial ingredient in proving the Brian\c{c}on-Skoda theorem in commutative algebra \cite{BS74, Huneke92, Li_BS} and an effective version of the Nullstellensatz in algebraic geometry \cite{EL99}. Moreover, a special case of Skoda's ideal generation also played key roles in Siu's works on the deformation invariance of plurigenera \cite{Siu98} and finite generation of the canonical ring \cite{Siu05}. The interaction between several complex variables, complex algebraic geometry and partial differential equations has been an attractive area for the researchers. For the sake of reader's convenience, we state a version of Skoda's $L^2$ division theorem as below.

\begin{theorem} \emph{(\cite{Skoda78}, Th\'eor\`eme 2).} \label{Skoda}
Let $(X,\omega)$ be an $n$-dimensional weakly pseudoconvex K\"ahler manifold with $\varphi\in\emph{Psh}(X)$, and $g:E\to Q$ be a surjective morphism of Hermitian holomorphic vector bundles with $r_E=\emph{rank}\,E$ and $r_Q=\emph{rank}\,Q$. Suppose that $E$ is Nakano semi-positive on $X$ and $L\to X$ is a Hermitian line bundle such that $$\sqrt{-1}\Theta(L)-\rho\sqrt{-1}\Theta(\det Q)\geq0$$ for $\rho=\min\{n,r_E-r_Q\}+\varepsilon$ and some $\varepsilon>0$.

Then, for every $f\in H^0(X,Q\otimes K_X\otimes L)$ satisfying
$$\int_{X}\langle\widetilde{gg^*}f,f\rangle\cdot(\det gg^*)^{-\rho-1}e^{-2\varphi}dV_\omega<+\infty,$$
there exists $h\in H^0(X,E\otimes K_X\otimes L)$ such that $f=g\cdot h$ and
$$\int_{X}|h|^2\cdot(\det gg^*)^{-\rho}e^{-2\varphi}dV_\omega \leq\frac{\rho}{\varepsilon}\cdot\int_{X}\langle\widetilde{gg^*}f,f\rangle\cdot(\det gg^*)^{-\rho-1}e^{-2\varphi}dV_\omega.$$
\end{theorem}

Due to Theorem \ref{Skoda}, if we consider the trivial bundles $E,\ Q$ and $L$ on a pseudoconvex domain, then by combining with the strong openness of multiplier ideal sheaves established by Guan-Zhou \cite{G-Z_open}, we can reformulate Theorem \ref{Skoda} in the language of multiplier ideals as follows (cf. also Remark \ref{Division_local_SGZ}):

\begin{theorem} \label{Skoda_MIS}
Let $X$ be an $n$-dimensional complex manifold with $\varphi\in\emph{QPsh}(X)$ a quasi-psh function and $\mathfrak{a}\subset\mathcal{O}_X$ a nonzero ideal sheaf with $r$ (local) generators. Then, it follows that
$$\mathscr{I}\big(\varphi+k\varphi_\mathfrak{a}\big)=\mathfrak{a}\cdot\mathscr{I}\big(\varphi+(k-1)\varphi_\mathfrak{a}\big),\ \forall k\geq\min\{n,r\},$$
where $\varphi_{\frak{a}}:=\frac{1}{2}\log(\sum_i|g_i|^2)$ and $(g_i)$ is any local system of generators of $\frak{a}$.
\end{theorem}

Motivated by the above reformulation of Theorem \ref{Skoda}, it is interesting for us to explore an analogue to Theorem \ref{Skoda_MIS} for Nadel-Lebesgue multiplier ideals in the singular setting. In order to achieve such a goal, a natural idea is to generalize Skoda's $L^2$ methods to the singular case, i.e., creating an appropriate $L^2$ theory for the $\overline\partial$-operator on singular complex spaces. However, as presented in \cite{Fornaess00, FG98}, it seems not to be possible to establish a general theory as in the smooth setting to solve the $\overline\partial$-equation with $L^2$ estimates on complex spaces with singularities; one can refer to \cite{FOV05, OV13, PS91, Ruppenthal14} for some partial results on the related topics.

On the other hand, we can also consider to apply Theorem \ref{Skoda} near the singularities under some reasonable assumptions on the positivity of curvatures. Fortunately, we could show that positivity of the Ricci curvature on the regular locus is in fact equivalent to the desired Skoda's ideal generation and $L^2$-solvability of the $\overline\partial$-equation. More precisely, we state our main result in the following:

\begin{theorem} \label{vanishing_regular}
Let $X$ be a (Hermitian) complex space of pure dimension $n$ with $x\in X$ a \emph{normal} point and $\pi:\widetilde X\to X$ a log resolution of the Jacobian ideal $\mathcal{J}ac_X$ of $X$. Then, the following statements are equivalent:
\begin{itemize}
  \item[(1)] For each quasi-psh function $\varphi$ near the point $x\in X$, we have
      $$R^q\pi_*\Big(\mathcal{O}_{\widetilde X}(\widehat K_{\widetilde X/X})\otimes\mathscr{I}(\varphi\circ\pi)\Big)=0,\ \forall q\geq1.$$

  \item[(2)] For each quasi-psh function $\varphi$ near the point $x\in X$, we have
      $$R^q\pi_*\Big(\mathcal{O}_{\widetilde X}(\widehat K_{\widetilde X/X})\otimes\mathscr{I}(\varphi\circ\pi)\Big)=0,\ \forall 1\leq q<n.$$

  \item[(3)] For some Stein neighborhood $\Omega\subset\subset X$ of the point $x$, there exists a K\"ahler metric $\omega$ on $\Omega$ such that the Ricci curvature $\emph{Ric}(\omega)\geq0$ on the regular locus $\Omega_\emph{reg}$ of $\Omega$.

  \item[(4)] For some Stein neighborhood $\Omega\subset\subset X$ of the point $x$, there exists a K\"ahler metric $\omega$ and a $\mathscr{C}^\infty$ differentiable real function $\psi$ on $\Omega$ such that $\emph{Ric}(\omega)+\sqrt{-1}\partial\overline{\partial}\psi\geq0$ on the regular locus $\Omega_\emph{reg}$ of $\Omega$.

  \item[(5)] For some Stein neighborhood $\Omega\subset\subset X$ of the point $x$, there exists a K\"ahler metric $\omega$ and a Hermitian line bundle $L$ on $\Omega$ such that for any smooth $\varphi\in\emph{SPsh}(\Omega)$ and $v\in L_{0,q}^2(\Omega_\emph{reg},L)$ satisfying $\overline\partial v=0$ and $$\int_{\Omega_\emph{reg}}\langle A_\varphi^{-1}v,v\rangle\,e^{-2\varphi}dV_\omega<+\infty$$ with the curvature operator $A_\varphi=[\sqrt{-1}\partial\overline{\partial}\varphi,\Lambda_\omega]$ on $\Omega_\emph{reg}$, we have $u\in L_{0,q-1}^2(\Omega_\emph{reg},L)$ such that $\overline\partial u=v$ and
      $$\int_{\Omega_\emph{reg}}|u|^2e^{-2\varphi}dV_\omega\leq\int_{\Omega_\emph{reg}}\langle A_\varphi^{-1}v,v\rangle\,e^{-2\varphi}dV_\omega.$$

  \item[(6)] For some Stein neighborhood $\Omega\subset\subset X$ of the point $x$, there exists a K\"ahler metric $\omega$ and a Hermitian line bundle
      $L$ on $\Omega$ such that for any smooth $\varphi\in\emph{SPsh}(\Omega)$ and $v\in L_{0,1}^2(\Omega_\emph{reg},L)$ satisfying $\overline\partial v=0$ and $$\int_{\Omega_\emph{reg}}\langle A_\varphi^{-1}v,v\rangle\,e^{-2\varphi}dV_\omega<+\infty,$$
      we have $u\in L^2(\Omega_\emph{reg},L)$ such that $\overline\partial u=v$ and
      $$\int_{\Omega_\emph{reg}}|u|^2e^{-2\varphi}dV_\omega\leq\int_{\Omega_\emph{reg}}\langle A_\varphi^{-1}v,v\rangle\,e^{-2\varphi}dV_\omega.$$

  \item[(7)] The Skoda's theorem holds for Nadel-Lebesgue multiplier ideals, i.e., for any nonzero ideal sheaf $\mathfrak{a}$ with $r$ generators and quasi-psh function $\varphi$ near the point $x\in X$, it holds that
      $$\mathscr{I}_\emph{NL}\big(\varphi+k\varphi_\mathfrak{a}\big)=\mathfrak{a}\cdot\mathscr{I}_\emph{NL}\big(\varphi+(k-1)\varphi_\mathfrak{a}\big),\ \forall k\geq\min\{n,r\}.$$

  \item[(8)] For any nonzero ideal sheaf $\mathfrak{a}$ near the point $x\in X$, it holds that
      $$\mathscr{I}_\emph{NL}\big(n\varphi_\mathfrak{a}\big)=\mathfrak{a}\cdot\mathscr{I}_\emph{NL}\big((n-1)\varphi_\mathfrak{a}\big).$$

  \item[(9)] For any nonzero ideal sheaf $\mathfrak{a}$ near the point $x\in X$, it holds that
      $$\mathscr{I}_\emph{NL}\big(n\varphi_\mathfrak{a}\big)\subset\mathfrak{a}.$$

  \item[(10)] The point $x\in X$ is a \emph{regular} point of $X$.
\end{itemize}
\end{theorem}

In the above result, the most interesting and amazing point is that it presents several characterizations of regular points by various statements involved Nadel-Lebesgue multiplier ideals, which look like almost irrelevant; e.g., (1, 2) in algebraic geometry, (3, 4) in differential geometry and (5, 6) in partial differential equations together with (7---9) in commutative algebra. The core idea of all arguments originates from the Skoda's ideal generation by the $L^2$ approaches in several complex variables.

\begin{remark}
Simple examples show that the assumption that $x$ is a normal point of $X$ cannot be removed in Theorem \ref{vanishing_regular}; in particular, any of the statements (1, 2, 7, 8) will not imply (10) in that case.
\end{remark}

As a straightforward consequence of Theorem \ref{vanishing_regular}, we have

\begin{corollary}
Any normal K\"ahler space with nonnegative Ricci curvature on the regular locus must be non-singular.
\end{corollary}

\subsection{K\"ahler-Einstein metrics on singular varieties}

Let $X$ be a normal $\mathbb{Q}$-Gorenstein K\"ahler space, that is, a normal K\"ahler space whose canonical class $K_X$ defines a $\mathbb{Q}$-line bundle on $X$. A K\"ahler current $\omega\in c_1(\pm K_X)$ is called a \emph{weak (or singular) K\"ahler-Einstein metric} on $X$ if $\omega$ has bounded local potentials and is a genuine K\"ahler-Einstein metric on the regular locus $X_\text{reg}$ of $X$ (cf. \cite{Berman16, BBEGZ, BG14, EGZ09}, etc.). A weak K\"ahler-Einstein metric $\omega$ on $X$ is called a K\"ahler-Einstein metric if $\omega$ is a K\"ahler metric on $X$, i.e., $\omega$ has smooth local potentials. For general expositions on the topic of K\"ahler-Einstein metrics one can refer to \cite{Aubin78, GZ_book17, Sz_book14, Xu21, Xu22, Yau78} and the
references therein. In particular, we state some recent results as follows.

\begin{theorem} \emph{(\cite{Berman16, BG14, EGZ09, Li22, LTW21, LTW22, LXZ22}, etc.).} \label{YTD_solution}
Let $X$ be a normal $\mathbb{Q}$-Gorenstein complex projective variety. Then:
\begin{itemize}
  \item[(1)] If $X$ is a $\mathbb{Q}$-Calabi-Yau variety with only log terminal singularities, then $X$ admits a weak K\"ahler-Einstein metric.

  \item[(2)] If $K_X$ is ample, then $X$ admits a weak K\"ahler-Einstein metric if and only if $X$ is $K$-stable.

  \item[(3)] If $-K_X$ is ample, then $X$ admits a weak K\"ahler-Einstein metric if and only if $X$ is $K$-polystable.
\end{itemize}
\end{theorem}

A basic and widely open problem in K\"ahler geometry/geometric analysis is understanding the geometric asymptotic behavior of the weak K\"ahler-Einstein metric near the singular locus $X_\text{sing}$ of $X$. In \cite{HS17}, the authors made a breakthrough with a very precise description for a class of Calabi-Yau varieties with smoothable isolated singularities, which are in further required to be isomorphic to a neighborhood of the vertex in a strongly regular Calabi-Yau cone; see also \cite{CS22, DFS21, FHJ21} for some recent progress in this direction. In more general situations, by using deep tools in the theory of degenerate complex Monge-Amp\`ere equations on singular complex spaces, the continuity of local potentials of weak K\"ahler-Einstein metrics is established for all $\mathbb{Q}$-Fano/Calabi-Yau varieties in \cite{BBEGZ, GGZ21}, but so far little is known for the higher order regularity in general and it is desirable to establish one for weak K\"ahler-Einstein potentials. However, relying on Theorem \ref{vanishing_regular}, we will see that too much regularity cannot be expected and in fact any weak K\"ahler-Einstein potential is at most $\mathscr{C}^\alpha\ (\alpha<2)$ differentiable near the singularities. In particular, we obtain the following

\begin{theorem} \label{KE_smooth}
Let $X$ be a normal $\mathbb{Q}$-Gorenstein K\"ahler space admitting a weak K\"ahler-Einstein metric $\omega$. Then, $\omega$ is smooth on $X$ if and only if $X$ is non-singular.
\end{theorem}

\section{Preliminaries}

Firstly, we introduce the notion of Nadel-Lebesgue multiplier ideal sheaf on any complex space of pure dimension and then present some useful facts used throughout this note.

\begin{definition}
Let $X$ be a complex space of pure dimension and $\varphi\in L_\text{loc}^1(X_\text{reg})$ with respect to the Lebesgue measure. Then, the complex space $X$ is said to be a \emph{Hermitian complex space} if there is a Hermitian metric $\omega$ on the regular part (may be disconnected) $X_\text{reg}$ of $X$ such that $\omega$ is locally the restriction of a Hermitian metric on some $\mathbb{C}^N$ for a local embedding of $X$. It follows from the differentiable partition of unity that every complex space is a Hermitian complex space as in the smooth case.

The complex space $X$ is called to be a \emph{K\"ahler space} if there is a Hermitian metric $\omega$ on $X$ such that $\omega$ is locally the restriction of a K\"ahler metric on some $\mathbb{C}^N$ for a local embedding of $X$. In particular, it admits smooth strictly psh functions as local potentials.

We say that the function $\varphi$ is \emph{quasi-plurisubharmonic} (quasi-psh for short) on $X$ if it is locally equal to the sum of a psh function and of a smooth function on $X$. The set of quasi-psh (resp. psh and strictly psh) functions on $X$ is denoted by $\text{QPsh}(X)$ (resp. $\text{Psh}(X)$ and $\text{SPsh}(X)$). A quasi-psh function $\varphi\in\text{QPsh}(X)$ will be said to have \emph{analytic singularities} on $X$ if $\varphi$ can be written locally as $$\varphi=\frac{c}{2}\log(|f_1|^2+\cdots+|f_{N_0}|^2)+O(1),$$
where $c\in\mathbb{R}_{\geq0}$ and $(f_i)$ are holomorphic functions.
\end{definition}

\begin{definition} \label{MIS_NL}
Let $(X,\omega)$ be a Hermitian complex space of pure dimension and $\varphi\in L_\text{loc}^1(X_\text{reg})$ with respect to the Lebesgue measure.

The \emph{Nadel-Lebesgue multiplier ideal sheaf} associated to $\varphi$ on $X$ is defined to be the $\mathcal{O}_X$-submodule $\mathscr{I}_\text{NL}(\varphi)\subset\mathscr{M}_X$ of germs of meromorphic functions $f\in\mathscr{M}_{X,x}$ such that $|f|^2e^{-2\varphi}$ is integrable with respect to the Lebesgue measure $dV_\omega$ near the point $x\in X$. One can check that $\mathscr{I}_\text{NL}(\varphi)$ is independent of the choice of Hermitian metric $\omega$ on $X$.

The \emph{log canonical threshold} (or \emph{complex singularity exponent}) $\text{LCT}_{x}(\varphi)$ of $\varphi$ at a point $x\in X$ is defined to be $$\text{LCT}_{x}(\varphi):=\sup\left\{c\ge0\ |\ \mathcal{O}_{X,x}\subset\mathscr{I}_\text{NL}(c\varphi)_x\right\}.$$
It is convenient to put $\text{LCT}_{x}(-\infty)=0$.

It is easy to see that $\mathscr{I}_\text{NL}(\varphi)\subset\mathcal{O}_X$ is an ideal sheaf when $X$ is a normal complex space and $\varphi$ is locally bounded from above on $X$. In addition, if $X$ is smooth and $\varphi\in\text{QPsh}(X)$, then $\mathscr{I}_\text{NL}(\varphi)$ is nothing but the usual multiplier ideal sheaf $\mathscr{I}(\varphi)$ introduced by Nadel (see \cite{De10}).
\end{definition}

\begin{remark}
Since the definition of Nadel-Lebesgue multiplier ideals is local, we can compute the multiplier ideals by choosing a special Hermitian metric $\omega$ for a local embedding of $X$. In particular, if $X$ is an $n$-dimensional complex subspace of some domain in $\mathbb{C}^N$, we can take Hermitian metric $\omega$ on $X$ to be the inherited standard K\"ahler metric from $\mathbb{C}^N$. Then, we have $dV_\omega=\frac{1}{n!}\upsilon^{n}|_{X_\text{reg}}$, where $\upsilon=\frac{\sqrt{-1}}{2}\sum\limits_{k=1}^{N}dz_k\wedge d\bar z_k$.
\end{remark}

For the sake of reader's convenience, we state a basic estimate related to local volume of an analytic subset as follows.

\begin{lemma} \emph{(\cite{G-L18}, Lemma 2.3).} \label{finiteness}
Let $X$ be a pure $n$-dimensional analytic subset through the origin $\textbf{\emph{0}}$ of some domain in $\mathbb{C}^N\ (N\geq2)$. Then, there is a Stein neighborhood $U\subset\subset\mathbb{C}^N$ of the origin $\textbf{\emph{0}}$ such that for any $0\leq\varepsilon<1$, we have
$$\int_{U\cap X}\frac{1}{(|z_1|^2+\cdots+|z_N|^2)^{n-1+\varepsilon}}dV_\omega<+\infty,$$
where $dV_\omega=\frac{1}{n!}\upsilon^n|_{X_\emph{reg}}\ \text{and}\ \upsilon=\frac{\sqrt{-1}}{2}\sum\limits_{k=1}^{N}dz_k{\wedge}d\bar{z}_k$.
\end{lemma}

Analogous to the Nadel-Ohsawa multiplier ideal sheaves introduced in \cite{Li_multiplier, Li_adjoint} (see also \cite{deFD14, EIM16} for the algebro-geometric counterpart), we state some related properties as follows.

\begin{proposition}  \label{MIS_NL_property}
$(1)$ Let $\pi:\widetilde X\to X$ be a log resolution of the Jacobian ideal $\mathcal{J}ac_X$ of $X$ and $\widehat K_{\widetilde X/X}$ be the Mather discrepancy divisor. Then, we have the image $$\emph{Im}\left(\pi^*\Omega_X^n\hookrightarrow\Omega_{\widetilde X}^n\right)=\mathcal{O}_{\widetilde X}(-\widehat K_{\widetilde X/X})\cdot\Omega_{\widetilde X}^n,$$
and
$$\mathscr{I}_\emph{NL}(\varphi)=\pi_*\left(\mathcal{O}_{\widetilde X}(\widehat K_{\widetilde X/X})\otimes\mathscr{I}(\varphi\circ\pi)\right).$$
Furthermore, we can deduce that $\mathscr{I}_\emph{NL}(\varphi+\log|\mathcal{J}ac_X|)=\mathscr{I}_\emph{NO}(\varphi)$, the Nadel-Ohsawa multiplier ideal sheaf associated to $\varphi$ on $X$.

$(2)$ When $X$ is normal and $\varphi$ has analytic singularities, $\mathscr{I}_\emph{NL}(\varphi)$ coincides with the Mather multiplier ideal sheaf defined in \cite{deFD14}.

$(3)$ For any $\varphi\in\emph{QPsh}(X)$, it follows that $\mathscr{I}_\emph{NL}(\varphi)\subset\mathscr{M}_X$ is a coherent fractional ideal sheaf and satisfies the strong openness, i.e., $\mathscr{I}_\emph{NL}(\varphi)=\bigcup\limits_{\varepsilon>0}\mathscr{I}_\emph{NL}\big((1+\varepsilon)\varphi\big)$.
\end{proposition}

For our proof of Theorem \ref{vanishing_regular}, we need the following $L^2$ estimates for the $\overline\partial$-equation and relative version of Grauert-Riemenschneider vanishing theorem for the higher direct images.

\begin{theorem} \emph{(cf. \cite{De10}, Theorem 5.2).} \label{dbar}
Let $(X,\omega)$ be an $n$-dimensional K\"ahler manifold, which contains a weakly pseudoconvex Zariski open subset. Let $L$ be a Hermitian line bundle on $X$ such that $\sqrt{-1}\Theta(L)+\emph{Ric}(\omega)>0$.

Then, for every smooth $\varphi\in\emph{Psh}(X)$ and $v\in L_{0,q}^2(X,L)$ satisfying $\overline\partial v=0$ and $$\int_{X}\langle A^{-1}v,v\rangle\,e^{-2\varphi}dV_\omega<+\infty$$
with the curvature operator $A=[\sqrt{-1}\Theta(L)+\emph{Ric}(\omega)+\sqrt{-1}\partial\overline{\partial}\varphi,\Lambda_\omega]$ on $X$,
there exists $u\in L_{0,q-1}^2(X,L)$ such that $\overline\partial u=v$ and
$$\int_{X}|u|^2e^{-2\varphi}dV_\omega\leq\int_{X}\langle A^{-1}v,v\rangle\,e^{-2\varphi}dV_\omega.$$
\end{theorem}

\begin{theorem} \emph{(\cite{DNWZ22}, Theorem 1.1).} \label{dbar_converse}
Let $(X,\omega)$ be an $n$-dimensional K\"ahler manifold which is a Zariski open subset of some Stein space $X^*$, and $L$ be a Hermitian line bundle on $X$.

If for any smooth $\varphi\in\emph{SPsh}(X^*)$ and $v\in L_{0,1}^2(X,L)$ satisfying $\overline\partial v=0$ and
$$\int_{X}\langle A_\varphi^{-1}v,v\rangle\,e^{-2\varphi}dV_\omega<+\infty$$
with the curvature operator $A_\varphi=[\sqrt{-1}\partial\overline{\partial}\varphi,\Lambda_\omega]$ on $X$,
there exists $u\in L^2(X,L)$ such that $\overline\partial u=v$ and
$$\int_{X}|u|^2e^{-2\varphi}dV_\omega\leq\int_{X}\langle A_\varphi^{-1}v,v\rangle\,e^{-2\varphi}dV_\omega,$$
then it follows that $L\otimes K_X^{-1}$ is Nakano semi-positive on $X$.
\end{theorem}

\begin{theorem} \emph{(\cite{Matsu_morphism}, Corollary 1.5).} \label{G-R_vanishing}
Let $\pi:X\to Y$ be a surjective proper (locally) K\"ahler morphism from a complex manifold $X$ to a complex space $Y$, and $(L,e^{-\varphi_L})$ be a (possibly singular) Hermitian line bundle on $X$ with semi-positive curvature. Then, the higher direct image sheaf
$$R^q\pi_*\Big(K_X\otimes L\otimes\mathscr{I}(\varphi_L)\Big)=0,$$
for every $q>\dim X-\dim Y$.
\end{theorem}

\begin{remark} \label{Re-K-Morphism}
Any log resolution $\pi:\widetilde X\to X$ of a coherent ideal sheaf $\mathscr{I}$ on a complex space $X$ is a locally K\"ahler (proper modification), which is locally a finite sequence of blow-ups with smooth centers. Besides, any finite holomorphic mapping between complex spaces is (locally) proper K\"ahler.
\end{remark}

In the remainder of this section, we recall some algebraic properties on the integral closure of ideals.

\begin{definition} (\cite{SH06}).
Let $R$ be a commutative ring and $I$ an ideal of $R$. An element $f\in R$ is said to be \emph{integrally dependent} on $I$ if it satisfies a relation
\[ f^d+a_1f^{d-1}+\cdots+a_d=0 \quad (a_k\in{I}^k, 1\le k\le d). \]

The set $\overline{I}$ consisting of all elements in $R$ which are integrally dependent on $I$ is called the \emph{integral closure} of $I$ in $R$. $I$ is called \emph{integrally closed} if $I=\overline{I}$. One can prove that $\overline{I}$ is an ideal of $R$, which is the smallest integrally closed ideal in $R$ containing $I$.
\end{definition}

\begin{definition} (\cite{SH06}).
Let $R$ be a commutative ring with identity and let $J\subset I$ be ideals in $R$. $J$ is said to be a \emph{reduction} of $I$ if there exists a nonnegative integer $n$ such that $I^{n+1}=JI^n$.

A reduction $J$ of $I$ is called \emph{minimal} if no ideal strictly contained in $J$ is a reduction of $I$. An ideal that has no reduction other than itself is called \emph{basic}.

One can prove that minimal reductions do exist in Noetherian local rings and an ideal which is a minimal reduction of a given ideal is necessarily basic. Moreover, if $R$ is a Noetherian ring, $J\subset I$ is a reduction of $I$ if and only if $\overline{J}=\overline{I}$.
\end{definition}

In the analytic setting, we have the following characterization on integral closure and reduction of ideals.

\begin{theorem} \emph{(cf. \cite{LJ-T}, Th\'eor\`eme 2.1).} \label{LJ-T_IC}
Let $X$ be a complex space and $Y\subset X$ be a proper closed complex subspace (may be non-reduced) defined by a coherent $\mathcal{O}_X$-ideal $\mathscr{I}$ with $x\in Y$ a point. Let $\mathscr{J}\subset\mathcal{O}_X$ be a coherent $\mathcal{O}_X$-ideal and $\mathcal{I}$ (resp. $\mathcal{J}$) be the germ of $\mathscr{I}$ (resp. $\mathscr{J}$) at $x$. Then, the following conditions are equivalent:
\begin{enumerate}
  \item $\mathcal{J}\subset\overline{\mathcal{I}}$.

  \item For every morphism $\pi:\widetilde X\to X$ satisfying: \emph{(i)} $\pi$ is a proper and surjective, \emph{(ii)} $\widetilde X$ is a normal complex space and \emph{(iii)} $\mathscr{I}\cdot\mathcal{O}_{\widetilde X}$ is an invertible $\mathcal{O}_X$-module, there exists an open neighborhood $U$ of $x$ in $X$ such that
      $$\mathscr{J}\cdot\mathcal{O}_{\widetilde X}|_{\pi^{-1}(U)}\subset\mathscr{I}\cdot\mathcal{O}_{\widetilde X}|_{\pi^{-1}(U)}.$$

  \item If $V$ is an open neighborhood of $x$ on which $\mathscr{I}$ and $\mathscr{J}$ are generated by their global sections, then for every system of generators $g_1,...,g_r\in\Gamma(V,\mathscr{I})$ and every $f\in\Gamma(V, \mathscr{J})$, one can find an open neighborhood $V'$ of $x$ and a constant $C>0$ such that
      $$|f(y)|\le C\cdot \sup_{k}|g_k(y)|,\ \forall y\in V'.$$
\end{enumerate}
\end{theorem}

\begin{remark} \label{IC_weaklyModi}
Let $X$ be a normal complex space and $\mathscr{I}\subset\mathcal{O}_X$ a coherent ideal sheaf. Let $\pi:\widetilde X\to X$ be any proper modification from a normal complex space $\widetilde X$ onto $X$ such that $\mathscr{I}\cdot\mathcal{O}_{\widetilde X}=\mathcal{O}_{\widetilde X}(-D)$ for some effective Cartier divisor $D$ on $\widetilde X$. Then, we have $\pi_*\mathcal{O}_{\widetilde X}(-D)=\overline{\mathscr{I}}$, the integral closure of $\mathscr{I}$ in $\mathcal{O}_X$.
\end{remark}

\begin{lemma} \emph{(cf. Example 9.6.19 in \cite{La04}; see also \cite{De10}, Lemma 11.16).} \label{Exis_reduction}
Let $X$ be a normal complex space of dimension $n$ and $\mathfrak{a}\subset\mathcal{O}_X$ a nonzero ideal. Then, there exists an open covering $\{U_\alpha\}_{\alpha\in\mathbb{N}}$ of $X$ such that $\mathfrak{a}|_{U_\alpha}$ has a reduction $\mathfrak{b}_\alpha$ generated by at most $n$ elements.
\end{lemma}

\section{Proofs of the main results}

\subsection{Proof of Theorem \ref{vanishing_regular}}

Since all of the statements are local, without loss of generality, we may assume that $X$ is an $n_{\geq2}$-dimensional normal (Hermitian) complex subspace of some domain in $\mathbb{C}^N$ with $\varphi\in\text{QPsh}(X)$ and $\mathfrak{a}=(g_1,\dots,g_r)\cdot\mathcal{O}_X$ an ideal sheaf generated by holomorphic functions $g_1,\dots,g_r$ on $X$.
Moreover, we may also assume that $\varphi$ is (locally) a strictly psh function on $X$ if necessary, by adding some smooth strictly psh function. It is easy to see that the implications $(1)\Longrightarrow(2),\ (3)\Longrightarrow(4),\ (5)\Longrightarrow(6)$ and $(7)\Longrightarrow(8)\Longrightarrow(9)$ are trivial; in particular, we will present a proof in the following order:

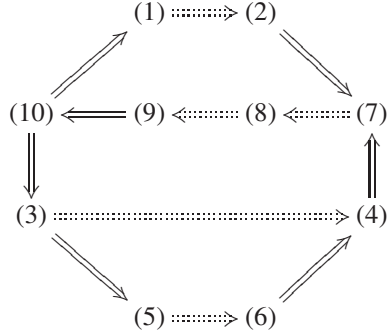
\[
\xymatrix
{
                                         & (1) \ar@{:>}[r]_{}  & (2) \ar@{=>}[dr]_{} &                   \\
    (10) \ar@{=>}[ur]_{}\ar@{=>}[d]_{}   & (9) \ar@{=>}[l]_{}  & (8) \ar@{:>}[l]_{}  & (7) \ar@{:>}[l]_{}\\
    (3)  \ar@{=>}[dr]^{}\ar@{:>}[rrr]^{} &                     &                     & (4) \ar@{=>}[u]_{}\\
                                         & (5) \ar@{:>}[r]_{}  & (6) \ar@{=>}[ur]_{} &
}
\]

\medskip

``$(2)\Longrightarrow(7)$''. By the definition of Nadel-Lebesgue multiplier ideal sheaf, it follows that
$$\mathfrak{a}\cdot\mathscr{I}_\text{NL}\big(\varphi+(k-1)\varphi_\mathfrak{a}\big)\subset\mathscr{I}_\text{NL}\big(\varphi+k\varphi_\mathfrak{a}\big),$$
and so it is sufficient to show the reverse inclusion.
\medskip

\textbf{Case (i).} When $r\leq n$.

Let $\pi:\widetilde X\to X$ be a common log resolution of $\mathcal{J}ac_X$ and $\mathfrak{a}$ such that $\mathfrak{a}\cdot\mathcal{O}_{\widetilde X}=\mathcal{O}_{\widetilde X}(-F)$ for some effective divisors $F$ on $\widetilde X$. Denote by
\begin{equation*}
\begin{split}
\mathscr{A}_m:=&\ \mathcal{O}_{\widetilde X}(\widehat K_{\widetilde X/X})\otimes\mathscr{I}(\varphi\circ\pi+m\varphi_\mathfrak{a}\circ\pi)\\
=&\ \mathcal{O}_{\widetilde X}(\widehat K_{\widetilde X/X})\otimes\mathscr{I}(\varphi\circ\pi)\otimes\mathcal{O}_{\widetilde X}(-mF)
\end{split}
\end{equation*}
for any $m\in\mathbb{N}$, and consider the Koszul complex determined by $g_1,\dots,g_r$:
$$0\to\Lambda^rV\otimes\mathcal{O}_{\widetilde X}(rF)\to\cdots\to\Lambda^2V\otimes\mathcal{O}_{\widetilde X}(2F)\to V\otimes\mathcal{O}_{\widetilde X}(F)\to\mathcal{O}_{\widetilde X}\to0,$$
where $V$ is the vector space spanned by $g_1,\dots,g_r$. Note that the Koszul complex is locally split and its syzygies are locally free, so twisting through by any coherent sheaf will preserve the exactness. Then, by twisting with $\mathscr{A}_k\ (k\geq r)$, we obtain the following long exact sequence
$$0\to\Lambda^rV\otimes\mathscr{A}_{k-r}\to\cdots\to\Lambda^2V\otimes\mathscr{A}_{k-2}\to V\otimes\mathscr{A}_{k-1}\to\mathscr{A}_k\to0.\eqno{(\star)}$$

On the other hand, for any $m\in\mathbb{N}$, by (2) we have the local vanishing of the higher direct images $R^q\pi_*\mathscr{A}_m=0\ (1\leq q<n)$. Note that
$$\mathscr{I}_\text{NL}\big(\varphi+m\varphi_\mathfrak{a}\big)=\pi_*\mathscr{A}_m$$
by the functoriality property with respect to direct images of sheaves by modifications, and then by taking direct
images of $(\star)$ we will deduce the following so-called exact Skoda complex (cf. \cite{La04}, p. 228):
$$0\to\Lambda^rV\otimes\mathscr{I}_\text{NL}\big(\varphi+(k-r)\varphi_\mathfrak{a}\big)\to\cdots\to V\otimes\mathscr{I}_\text{NL}\big(\varphi+(k-1)\varphi_\mathfrak{a}\big)\to\mathscr{I}_\text{NL}\big(\varphi+k\varphi_\mathfrak{a}\big)\to0.$$
In particular, the map $V\otimes\mathscr{I}_\text{NL}\big(\varphi+(k-1)\varphi_\mathfrak{a}\big)\to\mathscr{I}_\text{NL}\big(\varphi+k\varphi_\mathfrak{a}\big)$ is surjective, by which we can infer that $\mathscr{I}_\text{NL}\big(\varphi+k\varphi_\mathfrak{a}\big)\subset\mathfrak{a}\cdot\mathscr{I}_\text{NL}\big(\varphi+(k-1)\varphi_\mathfrak{a}\big)$ for any $k\geq r$.
\medskip

\textbf{Case (ii).} When $r>n$.

As the statement is local, then by Lemma \ref{Exis_reduction} we may assume that $\mathfrak{b}$ is a reduction of $\mathfrak{a}$ generated by $n$ elements $\widetilde g_1,...,\widetilde g_n$. Consider a common log resolution $\pi:\widetilde X\to X$ of $\mathcal{J}ac_X,\ \mathfrak{a}$ and $\mathfrak{b}$ such that $\mathfrak{a}\cdot\mathcal{O}_{\widetilde X}=\mathfrak{b}\cdot\mathcal{O}_{\widetilde X}=\mathcal{O}_{\widetilde X}(-F)$ for some effective divisors $F$ on $\widetilde X$. Then, by the same argument as above, we can deduce the following exact Skoda complex:
$$0\to\Lambda^nV\otimes\mathscr{I}_\text{NL}\big(\varphi+(k-n)\varphi_\mathfrak{a}\big)\to\cdots\to V\otimes\mathscr{I}_\text{NL}\big(\varphi+(k-1)\varphi_\mathfrak{a}\big)\to\mathscr{I}_\text{NL}\big(\varphi+k\varphi_\mathfrak{a}\big)\to0.$$
for any $k\geq n$, where $V$ is the vector space spanned by $\widetilde g_1,...,\widetilde g_n$. Therefore, it follows that $$\mathscr{I}_\text{NL}\big(\varphi+k\varphi_\mathfrak{a}\big)\subset\mathfrak{b}\cdot\mathscr{I}_\text{NL}\big(\varphi+(k-1)\varphi_\mathfrak{a}\big)\subset
\mathfrak{a}\cdot\mathscr{I}_\text{NL}\big(\varphi+(k-1)\varphi_\mathfrak{a}\big).$$
\medskip

``$(3)\Longrightarrow(5)$''. It follows from the assumption that we have a Stein neighborhood $\Omega\subset\subset X$ of the point $x$ with a K\"ahler metric $\omega$ such that $\text{Ric}(\omega)\geq0$ on $\Omega_\text{reg}$. Let $\varphi\in\text{SPsh}(\Omega)$ be any smooth strictly psh function on $\Omega$ and $L=\Omega\times\mathbb{C}$ be a trivial bundle equipped with the trivial Hermitian metric, which implies that $$\sqrt{-1}\Theta(L)+\text{Ric}(\omega)+\sqrt{-1}\partial\overline{\partial}\varphi\geq\sqrt{-1}\partial\overline{\partial}\varphi>0$$
on $\Omega_\text{reg}$.

Since $\Omega$ is a Stein space, we are able to choose a complex hypersurface $Z\subset\Omega$ which contains the singular locus $\Omega_\text{sing}$ of $\Omega$ such that $\Omega-Z\subset\Omega_\text{reg}$ is a Stein manifold. Then, by Theorem \ref{dbar} we obtain that, for any smooth $\varphi\in\text{SPsh}(\Omega)$ and $v\in L_{0,q}^2(\Omega_\text{reg},L)$ satisfying $\overline\partial v=0$ and
$$\int_{\Omega_\text{reg}}\langle A_\varphi^{-1}v,v\rangle\,e^{-2\varphi}dV_\omega<+\infty,$$
we can find $u\in L_{0,q-1}^2(\Omega_\text{reg},L)$ such that $\overline\partial u=v$ and
$$\int_{\Omega_\text{reg}}|u|^2e^{-2\varphi}dV_\omega\leq\int_{\Omega_\text{reg}}\langle A_\varphi^{-1}v,v\rangle\,e^{-2\varphi}dV_\omega.$$
\medskip

``$(6)\Longrightarrow(4)$''. As a straightforward application of Theorem \ref{dbar_converse} on $\Omega_\text{reg}$, it yields that $\sqrt{-1}\Theta(L)+\text{Ric}(\omega)\geq0$ on $\Omega_\text{reg}$. Let $\Omega'\subset\Omega$ be a small Stein neighborhood of the point $x$ such that the Hermitian line bundle $L$ has a smooth potential $\psi$ on $\Omega'$. Therefore, we deduce that $$\sqrt{-1}\partial\overline{\partial}\psi+\text{Ric}(\omega)=\sqrt{-1}\Theta(L)+\text{Ric}(\omega)\geq0$$ on $\Omega'_\text{reg}$.
\medskip

``$(4)\Longrightarrow(7)$''. Due to the definition and Lemma \ref{Exis_reduction}, it is sufficient to prove
$$\mathscr{I}_\text{NL}\big(\varphi+k\varphi_\mathfrak{a}\big)\subset\mathfrak{a}\cdot\mathscr{I}_\text{NL}\big(\varphi+(k-1)\varphi_\mathfrak{a}\big)$$
for the case $r\leq n$ near the point $x\in X$. Let $f\in\mathscr{I}_\text{NL}\big(\varphi+k\varphi_\mathfrak{a}\big)_x$ with $k\geq\min\{n,r\}=r$, then by the strong openness of multiplier ideals there exists small enough $\varepsilon>0$ such that $f\in\mathscr{I}_\text{NL}\big(\varphi+(k+\varepsilon)\varphi_\mathfrak{a}\big)_x$.

By the assumption of $(4)$, we let $\Omega\subset\subset X$ be a Stein neighborhood of the point $x$ with a K\"ahler metric $\omega$ and a smooth real function $\psi$ on $\Omega$ such that $\text{Ric}(\omega)+\sqrt{-1}\partial\overline{\partial}\psi\geq0$ on $\Omega_\text{reg}$. After shrinking $\Omega$ if necessary, we may assume that the function $\psi$ is bounded on $\Omega$ and $f$ is holomorphic on $\Omega$ such that
$$\int_{\Omega}|f|^{2}\cdot|g|^{-2(r+\varepsilon)}e^{-2\left(\varphi+(k-r)\varphi_\mathfrak{a}\right)}dV_\omega<+\infty.$$
In addition, we also choose a complex hypersurface $Z\subset\Omega$ which contains the singular locus $\Omega_\text{sing}$ of $\Omega$ and the common zero-set of holomorphic functions $g_1,...,g_r$ such that $\Omega':=\Omega-Z$ is a Stein manifold.

Let $E=\Omega'\times\mathbb{C}^r$ and $Q=\Omega'\times\mathbb{C}$ be the trivial bundles on $\Omega'$ and $L=K_{\Omega'}^{-1}$ be the anti-canonical line bundle with the induced metric twisted by a weight $e^{-\psi}$. The morphism $g:E\to Q$ determined by holomorphic functions $g_1,...,g_r$ is given by $$(h_1,...,h_r)\mapsto\sum\limits_{m=1}^rg_m\cdot h_m=g\cdot h.$$
Note that $\widetilde{gg^*}=\text{Id}_Q$ when $\text{rank}\,Q=1$, and on $\Omega'$ we have
$$\sqrt{-1}\Theta(L)-(r-1+\varepsilon)\sqrt{-1}\Theta(\det Q)=\text{Ric}(\omega)+\sqrt{-1}\partial\overline{\partial}\psi\geq0.$$
Thus, we can apply Theorem \ref{Skoda} on $\Omega'$ and then obtain an $r$-tuple $(h_1,...,h_r)$ of holomorphic functions on $\Omega'$ such that $f=g\cdot h$ on $\Omega'$ and
$$\int_{\Omega'}|h|^{2}\cdot|g|^{-2(r-1+\varepsilon)}e^{-2\left(\varphi+(k-r)\varphi_\mathfrak{a}\right)}dV_\omega =\int_{\Omega'}|h|^{2}e^{-2\left(\varphi+(k-1+\varepsilon)\varphi_\mathfrak{a}\right)}dV_\omega<+\infty.$$
We can now extend every $h_m$ to be a holomorphic function on $\Omega$ from the $L^2$ estimate above and normality of $X$, which implies that
$$\mathscr{I}_\text{NL}\big(\varphi+k\varphi_\mathfrak{a}\big)\subset\mathfrak{a}\cdot\mathscr{I}_\text{NL} \big(\varphi+(k-1+\varepsilon)\varphi_\mathfrak{a}\big)\subset\mathfrak{a}\cdot\mathscr{I}_\text{NL}\big(\varphi+(k-1)\varphi_\mathfrak{a}\big)$$
on $\Omega$; we finish the argument.
\medskip

``$(9)\Longrightarrow(10)$''. By the assumption, we have $\mathscr{I}_\text{NL}\big(n\varphi_\mathfrak{a}\big)\subset\mathfrak{a}$. Suppose that $x\in X$ is a singular point. Then, by the local parametrization for analytic sets, we can find a local coordinate system $(z';z'')=(z_1,...,z_n;z_{n+1},...,z_N)$ near $x$ such that for some constant $C>0$, we have $|z''|\leq C\cdot|z'|$ for any point $z\in X$ near $x$.

Let $\mathfrak{a}\subset\mathcal{O}_X$ be the ideal sheaf generated by holomorphic functions $\widehat{z}_1,...,\widehat{z}_{n}\in\mathcal{O}_X$ (shrinking $X$ if necessary), where $\widehat{z}_k$ are the residue classes of $z_k$ in $\mathcal{O}_X$. From the non-smoothness of $X$ at the point $x$, we deduce that the embedding dimension $\dim_{\mathbb{C}}(\mathfrak{m}_{X,x}/\mathfrak{m}^2_{X,x})\geq n+1$ of $X$ at $x$, which implies that there exists $k_0\ (n+1\leq k_0\leq N)$ such that $\widehat{z}_{k_0}\not\in\mathfrak{a}$.

On the other hand, after shrinking $X$ again, it follows that
$$\int_X\frac{|{z}_{k_0}|^2}{|z'|^{2n}}dV_\omega\leq C^2(1+C^2)^{n-1}\cdot\int_X|z|^{-2(n-1)}dV_\omega<+\infty,$$
where the finiteness of the integration follows from Lemma \ref{finiteness}. Then, we infer that $\widehat{z}_{k_0}\in\mathscr{I}_\text{NL}\big(n\varphi_\mathfrak{a}\big)$, but $\widehat{z}_{k_0}\not\in\mathfrak{a}$, which contradicts to the assumption $\mathscr{I}_\text{NL}\big(n\varphi_\mathfrak{a}\big)\subset\mathfrak{a}$.
Thus, we obtain that $x\in X$ is a regular point.
\medskip

``$(10)\Longrightarrow(1)$''. It is a straightforward consequence of Theorem \ref{G-R_vanishing}.
\medskip

``$(10)\Longrightarrow(3)$''. Since $x$ is a regular point of $X$, after choosing an appropriate coordinate neighborhood of $x$, we may assume that $\Omega\ni x$ is a Stein domain in $\mathbb{C}^n$. Therefore, we can take $\omega=\frac{\sqrt{-1}}{2}\sum\limits_{k=1}^{n}dz_k\wedge d\bar z_k$ to be the standard Euclidean metric on $\mathbb{C}^n$ and then we have $\text{Ric}(\omega)=0$ on $\Omega$; the proof of Theorem \ref{vanishing_regular} is concluded.
\hfill $\Box$

\begin{remark}
In addition, we can deduce from the proof of Theorem \ref{vanishing_regular} that

(i) if $(1)$ or $(2)$ holds for each quasi-psh function $\varphi$ with analytic singularities, then $x\in X$ is a regular point;

(ii) both of the statements $(3)$ and $(4)$ could be respectively modified to be $\text{Ric}(\omega)\geq0$ and $\text{Ric}(\omega)+\sqrt{-1}\partial\overline{\partial}\psi\geq0$ on a Zariski open subset of $\Omega$ contained in $\Omega_\text{reg}$.
\end{remark}

\subsection{Proof of Theorem \ref{KE_smooth}}

It is sufficient to prove the necessity.

Let $x\in X$ be any point. Since $\omega$ is a smooth K\"ahler metric on $X$, then $\omega$ has smooth local potentials, i.e., there exists a Stein neighborhood $\Omega\subset X$ of $x$ and a smooth strictly psh functions $\psi$ on $\Omega$ such that $\omega=\sqrt{-1}\partial\overline{\partial}\psi$ on $\Omega_\text{reg}$, which implies that $\text{Ric}(\omega)+\sqrt{-1}\partial\overline{\partial}\psi\geq0$ on $\Omega_\text{reg}$ whenever $\text{Ric}(\omega)=\pm\omega,\,0$. Thus, it follows from $(4)$ in Theorem \ref{vanishing_regular} that $x\in X$ is a regular point.
\hfill $\Box$
\smallskip

\begin{remark}
The same arguments as in the proof of Theorem \ref{vanishing_regular} and \ref{KE_smooth} also imply that each local potential of weak K\"ahler-Einstein metric $\omega$ is $\mathscr{C}^2$ differentiable on $X$ if and only if $X$ is non-singular, and that there exists no \emph{singular} normal K\"ahler space such that the K\"ahler metric is K\"ahler-Einstein on the regular locus.

In fact, our method is still available when the weak K\"ahler-Einstein metric is (locally) equivalent to the standard induced K\"ahler metric by restriction near the singularities; for instance, when the weak K\"ahler-Einstein metric is of locally bounded coefficients.
\end{remark}

\appendix
   \renewcommand{\appendixname}{Appendix~\Alph{section}}

\section{Uniform bounds of powers associated to an $L^2$ division problem} \label{A}

The ideal membership is an important object to study in commutative algebra, algebraic geometry and several complex variables, e.g., the famous Hilbert's Nullstellensatz and Brian\c{c}on-Skoda theorem and so on. In this part, we are mainly interested in the uniform bounds of powers associated to an $L^2$ division problem, a kind of special ideal membership. Let $X$ be a Stein manifold of dimension $n$ and $\mathfrak{a}=(g_1,\dots,g_r)\cdot\mathcal{O}_X$ an ideal sheaf generated by holomorphic functions $g_1,\dots,g_r$ on $X$. In general, the division problem states that, given a positive integer $k\in\mathbb{N}$ and holomorphic function $f$ on $X$, we wish to determine when $f$ is generated by holomorphic functions $g_1,\dots,g_r$; more precisely, when we can find holomorphic functions $h_1,\dots,h_r\in\mathfrak{a}^{k-1}$ on $X$ such that $$f=\sum\limits_{m=1}^rg_m\cdot h_m.$$
Thanks to the Oka-Cartan's theory on Stein manifolds, the division problem is solvable if and only if $f\in\mathfrak{a}^k$.

Note that the condition $f\in\mathfrak{a}^k$ is purely algebraic, and so it is natural to ask whether we could find an analytic condition to replace the algebraic one. It is easy to see that $f\in\mathfrak{a}^k$ implies that $|f|e^{-\varphi_k}$ is locally bounded on $X$, or $L_{\text{loc}}^2$ more generally; where $\varphi_k:=k\log|g|$ and $|g|^2:=|g_1|^2+\cdots+|g_r|^2$. On the other hand, local boundedness of $|f|e^{-\varphi_k}$ is equivalent to the fact that $f\in\overline{\mathfrak{a}^k}$, the integral closure of $\mathfrak{a}^k$ in $\mathcal{O}_X$ (see Theorem \ref{LJ-T_IC}). Thus, it is an interesting question whether we could establish solvability of an $L^2$ analogue of the division problem.

Let $\varphi\in\text{Psh}(X)$ be a psh function on $X$ and denote by
$$A_\text{loc}^2(X,\varphi):=\left\{f\in\mathcal{O}_X(X)\ \big|\ |f|^2e^{-2\varphi}\ \text{is locally integrable on $X$}\right\}.$$
Then, we raise the following $L^2$ division problem:

\begin{question} \label{Q_Division}
Let $X$ be an $n$-dimensional Stein manifold with a psh function $\varphi\in\emph{Psh}(X)$, and $\mathfrak{a}=(g_1,\dots,g_r)\cdot\mathcal{O}_X$ an ideal sheaf generated by holomorphic functions $g_1,\dots,g_r$ on $X$. Given positive integer $k\in\mathbb{N}$ and $f\in A_\emph{loc}^2(X,\varphi+\varphi_k)$, are there holomorphic functions $h_1,\dots,h_r\in A_\emph{loc}^2(X,\varphi+\varphi_{k-1})$ such that $$f=\sum\limits_{m=1}^rg_m\cdot h_m?$$
\end{question}

\subsection{A solution to Question \ref{Q_Division}}

Unfortunately, the answer of Question \ref{Q_Division} is negative for general $k$ (see Example \ref{Division_neg}). Motivated by the Skoda's $L^2$ division theorem (cf. Theorem \ref{Skoda}), it seems to be reasonable to find a uniform integer $k_0$, depending only on $n$, such that Question \ref{Q_Division} is solvable for any $k\geq k_0$. The goal of this part is to present an optimal uniform lower bounds of powers associated to Question \ref{Q_Division}. In particular, we will establish the following

\begin{theorem} \label{Division_local}
There exists a uniform integer $k_0=\min\{n,r\}$ such that the solution to Question \ref{Q_Division} is positive for any $k\geq k_0$. In further, the uniform lower bound $k_0=\min\{n,r\}$ is optimal.
\end{theorem}

In fact, the optimality of uniform integer $k_0=\min\{n,r\}$ is straightforward by the following:

\begin{example} \label{Division_neg}
Let $B_n(\textbf{\emph{0}})$ be the unit ball centered at the origin $\textbf{\emph{0}}=(\textbf{\emph{0}}',\textbf{\emph{0}}'')$ in $\mathbb{C}^r\times\mathbb{C}^{n-r}\ (1\leq r\leq n)$ and take $g_1=z_1,...,g_r=z_r,f\equiv1,\varphi\equiv0$ on $B_n(\textbf{\emph{0}})$. Then, for every $k<k_0=r$, the answer of Question \ref{Q_Division} is negative.

Indeed, by the fact that the log canonical threshold $\emph{LCT}_{(\textbf{\emph{0}}',z'')}(\varphi_k)=\frac{r}{k}>1$ of $\varphi_k$ at any point $(\textbf{\emph{0}}',z'')$, one can derive that $f\in A_\emph{loc}^2(B_n(\textbf{\emph{0}}),\varphi+\varphi_k)$. Then, we infer from the fact that $f$ has no zeros in $B_n(\textbf{\emph{0}})$ that there exist no holomorphic functions $h_1,\dots,h_r$ on $B_n(\textbf{\emph{0}})$ such that $f=\sum\limits_{m=1}^rg_m\cdot h_m$.
\end{example}

\noindent{\bf{Proof of Theorem} \ref{Division_local}.} It follows from the local vanishing (Theorem \ref{G-R_vanishing}) and the arguments as in the proof of Theorem \ref{vanishing_regular} that for any $k\geq\min\{n,r\}$, we have
$$\mathscr{I}\big(\varphi+\varphi_k\big)=\mathfrak{a}\cdot\mathscr{I}\big(\varphi+\varphi_{k-1}\big).$$

Let $$\tau:\mathscr{I}\big(\varphi+\varphi_{k-1}\big)^{\oplus r}\longrightarrow\mathscr{I}\big(\varphi+\varphi_k\big)$$ be the sheaf homomorphism defined by $$\tau(h_{1,x},\dots,h_{r,x})=\sum\limits_{m=1}^rg_m\cdot h_{m,x}$$ for any germs $h_{m,x}\in\mathscr{I}\big(\varphi+\varphi_{k-1}\big)_x$. Then, we have an exact sequence of sheaves
$$\mathscr{I}\big(\varphi+\varphi_{k-1}\big)^{\oplus r}\stackrel{\tau}{\longrightarrow}\mathscr{I}\big(\varphi+\varphi_k\big)\longrightarrow0.$$
It follows from the Oka-Cartan theory on Stein manifolds that the induced sequence of sections
$$\Gamma\Big(X,\mathscr{I}\big(\varphi+\varphi_{k-1}\big)^{\oplus r}\Big)\stackrel{\tau^*}{\longrightarrow}\Gamma\Big(X,\mathscr{I}\big(\varphi+\varphi_k\big)\Big)\longrightarrow0$$
is also exact, which implies that any section $f\in\Gamma\Big(X,\mathscr{I}\big(\varphi+\varphi_k\big)\Big)$ can be written as the image $f=\sum\limits_{m=1}^rg_m\cdot h_m$ for some sections $h_m\in\Gamma\Big(X,\mathscr{I}\big(\varphi+\varphi_{k-1}\big)\Big)$.
\hfill $\Box$

\begin{remark} \emph{(An alternative argument on Theorem \ref{Division_local}).} \label{Division_local_SGZ}
In fact, we could also give another argument on the proof of Theorem \ref{Division_local} depending on the strong openness of multiplier ideals established by Guan-Zhou \cite{G-Z_open} and the Skoda's $L^2$ division theorem for holomorphic functions (see Theorem \ref{Skoda}).

Since the statement is local, it follows from Lemma \ref{Exis_reduction} that it is sufficient to prove $\mathscr{I}\big(\varphi+\varphi_k\big)\subset\mathfrak{a}\cdot\mathscr{I}\big(\varphi+\varphi_{k-1}\big)$ for the case $r\leq n$. Given $f\in\Gamma\Big(X,\mathscr{I}\big(\varphi+\varphi_k\big)\Big)$, after shrinking $X$, we may assume that $X$ is the unit ball in $\mathbb{C}^n$ and $$\int_X|f|^2e^{-2(\varphi+\varphi_k)}d\lambda_n=\int_X|f|^2\cdot|g|^{-2k}e^{-2\varphi}d\lambda_n<+\infty.$$
Then, for each $k\geq r$, by the strong openness of multiplier ideals there exists sufficiently small $\varepsilon>0$ such that
$$\int_X|f|^2e^{-2(\varphi+(1+\varepsilon)\varphi_k)}d\lambda_n=\int_X|f|^2\cdot|g|^{-2(1+\varepsilon)k}e^{-2\varphi}d\lambda_n<+\infty,$$
shrinking $X$ if necessary. Finally, combining with Theorem \ref{Skoda}, we deduce the desired result.
\end{remark}

\subsection{A global $L^2$ version of Question \ref{Q_Division}}

Let $(X,\omega)$ be an $n$-dimensional Stein manifold with a K\"ahler form $\omega$. Let $\varphi\in\text{Psh}(X)$ and $\mathscr{I}=(g_1,\dots,g_r)\cdot\mathcal{O}_X$ an ideal sheaf generated by holomorphic functions $g_1,\dots,g_r$ on $X$. Denote by $$A^2(X,\varphi):=\left\{f\in\mathcal{O}_X(X)\ \bigg|\ \int_X|f|^2e^{-2\varphi}dV_\omega<+\infty\right\}.$$
Then, we have the following global analogue of Question \ref{Q_Division}:

\begin{question} \label{Q_DivisionG}
Can we find a uniform integer $k_0$ such that for each $k\geq k_0$ and $f\in A^2(X,\varphi+\varphi_k)$, there exist $h_1,\dots,h_r\in A^2(X,\varphi+\varphi_{k-1})$ satisfying $$f=\sum\limits_{m=1}^rg_m\cdot h_m?$$
\end{question}

As an immediate consequence of Theorem \ref{Skoda}, we obtain the following

\begin{theorem} \label{Division_global}
Let $X$ be a pseudoconvex domain in $\mathbb{C}^n$. Then, there exists a uniform integer $k_0=\min\{n+2,r+1\}$ such that the solution to Question \ref{Q_DivisionG} is positive.
\end{theorem}

\begin{remark}
$(1)$ More generally, Theorem \ref{Division_global} also holds for any complete K\"ahler domain in $\mathbb{C}^n$ with smooth psh function $\varphi\in\text{Psh}(X)$.

$(2)$ In this case, combining with the Example \ref{Division_neg}, it follows that the optimal uniform lower bound $k_0$ is at least $\min\{n,r\}$, and at most $\min\{n+2,r+1\}$.
\end{remark}

\end{document}